%% file: PflPosTanGGCDS.tex
\documentclass[11]{amsart}

\usepackage{amsmath,amsfonts,amsxtra}
\usepackage{amssymb}
\let\amslrcorner\lrcorner
\usepackage{hyperref}
\usepackage{mathtools}
\usepackage{mathrsfs}
\usepackage{enumerate}
\usepackage{comment}
\usepackage{mathabx}

\usepackage[width=5.8in, height=8.5in, bottom=1.3in, margin=1.2in, centering]{geometry}
\numberwithin{equation}{section}

\theoremstyle{plain}
        \newtheorem{theorem}{Theorem}[section]
        \newtheorem{lemma}[theorem]{Lemma}

        \theoremstyle{definition}
        \newtheorem{definition}[theorem]{Definition}
        \newtheorem{remark}[theorem]{Remark}
        \newtheorem{example}{Example}[section]

\newcommand{\N}{{\mathbb N}}

\newcommand{\C}{{\mathbb C}}
\newcommand{\R}{{\mathbb R}}

\newcommand{\frakg}{\mathfrak{g}}
\newcommand{\frakh}{\mathfrak{h}}
\newcommand{\frakm}{\mathfrak{m}}

\newcommand{\calC}{\mathcal{C}}
\newcommand{\calO}{\mathcal{O}}

\newcommand{\scrA}{\mathscr{A}}
\newcommand{\scrC}{\mathscr{C}}

\newcommand{\scrI}{\mathscr{I}}
\newcommand{\scrJ}{\mathscr{J}}

\newcommand{\scrB}{\mathscr{B}}

\newcommand{\scrL}{\mathscr{L}}
\newcommand{\scrO}{\mathscr{O}}
\newcommand{\scrQ}{\mathscr{Q}}
\newcommand{\scrS}{\mathscr{S}}
\newcommand{\Zeta}{\mathrm{Z}}

\newcommand{\wtildeU}{\widetilde{U}}

\newcommand{\whatU}{\widehat{U}}

\newcommand{\spec}{\operatorname{spec}}
\newcommand{\sphere}{\mathrm{S}}
\newcommand{\inertianull}{{\Lambda_0}}
\newcommand{\sttimes}{{\hspace{0.1em}{_s}\hspace{-.1em}\times_{\hspace{-.05em}t}\hspace{0.1em}}}
\newcommand*{\longhookrightarrow}{\ensuremath{\lhook\joinrel\relbar\joinrel\rightarrow}}

\title[Grauert--Grothendieck complex on differentiable spaces]{The Grauert--Grothendieck complex on differentiable spaces and a sheaf complex of Brylinski}

\author{Markus J.~Pflaum, Hessel Posthuma,~\textrm{and} Xiang Tang}
\address{\newline Markus J. Pflaum, 
         {\rm Department of Mathematics, University of Colorado UCB 395,
         Boulder CO 80309, USA} \newline\indent
         {\rm email:} {\tt markus.pflaum@colorado.edu}\newline
        Hessel Posthuma,
         {\rm Korteweg-de Vries Institute for Mathematics, University of Amsterdam, The Netherlands} 
         \newline\indent
         {\rm email:} {\tt H.B.Posthuma@uva.nl}\newline
        Xiang Tang,  
         {\rm  Department of Mathematics, Washington University, St.~Louis, USA}
         \newline\indent
         {\rm email:} {\tt xtang@math.wustl.edu} 
      }

\begin{document}
\maketitle
\dedicatory{\begin{center}
Dedicated to Henri Laufer on the occasion of his 70th birthday.
\end{center}
}

\begin{abstract}
We use the Grauert--Grothendieck complex on differentiable 
spaces to study basic relative forms on the inertia space of a compact
Lie group action on a manifold. We prove that the sheaf complex of basic 
relative 
forms on the inertia space is a fine resolution of Bryliski's sheaf of 
functions on the inertia space.
\end{abstract}

\input{intro}
\input{forms}
\input{basic-relative-forms}

\appendix
\input{appendix}
\bibliographystyle{amsalpha}    
\bibliography{PflPosTanGGCDS}

\end{document}

%% file: intro.tex
\section{Introduction}

In his paper \cite{BryCHET}, {\sc Jean-Luc Brylinski} studies the cyclic homology theory of 
the smooth crossed-product algebra $\scrA := \scrC^\infty_\textup{c} (G \times M)$ associated to a manifold 
$M$ which carries a smooth action of a Lie group $G$. 
The crossed-product algebra carries the convolution product $\ast$ defined by 
\begin{equation}
  \label{eq:convolution-product}
  \big( u \ast v \big) ( g ,x ) := \int_M  u(gh,h^{-1}x) \cdot v(h^{-1} , x)  \, d\mu (h) 
  \quad \text{for $u,v \in \scrA$, $(g,x) \in G \times M$.} 
\end{equation}
The symbol $\mu$ hereby denotes a fixed left-invariant Haar measure on $G$.
{\sc Brylinski} asserts in his article that the Hochschild 
homology $HH_k (\scrA)$ of the (topological) algebra $\scrA$ coincides 
naturally with the space $\Omega^k_{\textup{b.r.}} (M/G)$ of so-called \emph{basic  relative 
$k$-forms}. The sheaf complex of basic relative forms will be constructed in 
Section \ref{sec:basic-relative-forms}. Despite it being a sheaf complex over the  orbit space 
$M/G$,  basic relative forms are defined as forms on the so-called 
\emph{loop space}
\begin{equation}
  \label{eq:loop-space}
  \inertianull (G \ltimes M) := \bigcup_{g \in G} \{ g \} \times M^g \subset G \times M
\end{equation}
which essentially consists of a disjoint union of the fixed point manifolds 
$M^g$, $g\in G$ of the smooth action $G \times M \to M$ together with an 
appropriate topology on it. In Section \ref{sec:basic-relative-forms} we will also
see that the loop space carries even the structure of a differentiable stratified space.
{\sc Brylinski} has also claimed in \cite{BryCHET} that the
sheaf complex of basic relative forms is a resolution of a certain sheaf $\scrB$ on the orbit space $M/G$.   
The section spaces of that sheaf over $O \subset M/G$ open are given by
\begin{equation}
  \label{eq:brylinski-sheaf-section-spaces}
  \scrB(O) := \big\{ f \in \scrC^\infty_{\inertianull} (s^{-1}_{|\inertianull} \pi^{-1} (O) ) \mid f 
  \text{ is $G$-invariant and $f(g,-)$ locally constant for all $g\in G$}  \big\} .
\end{equation}
Hereby, the map $s_{|\inertianull} : \inertianull (G\ltimes M) \to M$ is given by 
$s_{|\inertianull}(g,p)=p$, $\pi : M \to G/M$ is the orbit projection, 
the sheaf $\scrC^\infty_{\inertianull}$ is defined as the sheaf of continuous functions
on the loop space which are locally restrictions of smooth functions on $G\times M$, 
the $G$-action on $\inertianull (G \ltimes M)$  is given by the diagonal action with 
conjugation in the first coordinate, and the function $f(g,-)$ for 
$f\in  \scrC^\infty (\inertianull (G\ltimes M))$ and $g\in G$ is the map 
$M^g \to \R$, $p \to f(g,p)$. Since the sheaf $\scrB$ has been defined first in 
\cite{BryCHET}, we call it \emph{Brylinskis's sheaf}.

For the proof that $HH_k (\scrA)$ is isomorphic to $\Omega^k_\textup{b.r.} (M/G)$, {\sc Brylinski} refers 
to the unpublished paper \cite{BryAAGAH}, a proof of the second claim is missing.

The purpose of these notes is to shed some light onto Brylinski's
sheaf $\scrB$ and the sheaf complex of basic relative forms. We interpret the latter as a certain 
subcomplex of the Grauert-Grothendieck complex of differential forms 
on the loop space of the $G$-manifold $M$. Then, in Theorem \ref{thm:ResBrylinskiSheaf},we show that it is what 
it is claimed to be: 
an acyclic complex of fine sheaves whose cohomology in degree zero coincides 
with Brylinski's sheaf $\scrB$.

Let us mention that the cyclic homology theory of crossed product algebras has been studied 
also by {\sc Block--Getzler} \cite{BloGetECHEDF} and {\sc Nistor} \cite{NisCCCPAG}. Moreover, {\sc Farsi--Pflaum--Seaton} have described in \cite{FarsiPflaumSeaton,FarPflSeaDSGDRTIS} the stratification theory of the loop space and proved 
a de Rham Theorem for it.
\vspace{2mm}

\noindent
{\bf Acknowledgments:}
M.P.~acknowledges financial support by the Simons Foundation under 
Collaboration Grant nr.~359389.  
He also thanks the Max-Planck-Institut for Mathematics of the Sciences in 
Leipzig, Germany, and the Tsinghua Sanya International Mathematics Forum, 
China, for hospitality and support.

X.T.~acknowledges funding by the NSF under award DMS 1363250. He also 
would like to thank the University of Amsterdam, Korteweg-de Vries 
Institute for Mathematics, and the Max-Planck Institute at Leipzig for hosting 
his visits, where part of the work has been done.

%% file: forms.tex
\section{Differential forms}
\label{SectionDiffForms}
In this section, we will define the Grauert--Grothendieck complex of 
differential forms on a differentiable space $(X,\scrO)$ which we  
allow to be non-reduced. To this end consider first an affine open subset 
of $X$ and let $\iota: U \hookrightarrow \wtildeU\subset \R^n$ be a singular chart. 
Denote by $\scrI_\iota$ the ideal sheaf in $\scrC^\infty_{\wtildeU}$ such that 
$ \scrO_{|U} \cong \iota^{-1} (\scrO_\iota)$, where $\scrO_\iota$ is the restriction of the sheaf 
$ \scrC^\infty_{\wtildeU} / \scrI_\iota$ to $\iota (U)$. 
We then define $\Omega^0_\iota$ as the sheaf $\scrO_\iota$ and, for $k \in \N^*$,  
the sheaf $\Omega^k_\iota$ as  the inverse image sheaf 
\begin{equation}
\label{eq:def-diff-k-forms}
  \Omega^k_\iota  := \iota^{-1}\big( \Omega^k_{\wtildeU} / \scrI_\iota \Omega^k_{\wtildeU}  + 
  d \scrI_\iota \wedge \Omega^{k-1}_{\wtildeU} \big) = 
  \iota^*\big( \Omega^k_{\wtildeU} / \scrI_\iota \Omega^k_{\wtildeU}  + 
  d \scrI_\iota \wedge \Omega^{k-1}_{\wtildeU} \big) \: .
\end{equation}
The latter equality holds true because 
$\scrO_{|U} \cong \iota^{-1} (\Omega^0_\iota) $, and because the sheaf 
\[
  \Omega^k_{\wtildeU} / \scrI_\iota \Omega^k_{\wtildeU}  + d \scrI_\iota \wedge \Omega^{k-1}_{\wtildeU}
\]
is an $\scrO_\iota$-module. By construction, the exterior differential $d$ 
descends to sheaf morphisms $d: \Omega^k_\iota \to  \Omega^{k-1}_\iota $ so that we obtain
a complex of sheaves $\big( \Omega^\bullet_\iota , d \big)$ over $U$. In the following we will
show that the sheaf complexes $\big( \Omega^\bullet_\iota , d \big)$ glue to a sheaf complex
$\big( \Omega^\bullet , d \big)$ on $X$, when $\iota$ runs through a singular atlas of $(X,\scrO)$. 
To construct the gluing maps we will proceed by proving a sequence of lematas.

\begin{lemma}
\label{lem:step1}
  Let $\iota: U \hookrightarrow \wtildeU\subset \R^n$ and $\kappa: U \hookrightarrow \whatU \subset \R^m$ 
  be two singular charts of $(X,\scrO)$ for which there exists a smooth embedding 
  $H:\whatU \hookrightarrow \wtildeU$ such that $H(\whatU)$ is closed in 
  $\wtildeU$ and such that the pullback $H^* : \scrC^\infty_{\wtildeU} \to \scrC^\infty_{\widehat{U}}$
  induces an isomorphism of locally ringed spaces 
  $\big( H_{|\kappa (U)},H^*\big) : (\kappa (U) , \scrO_\kappa ) \to  
  (\iota (U) , \scrO_\iota  ) $. 
  %
  Then there exists a unique isomorphism of sheaf complexes 
  $ \eta_{\kappa,\iota} : \Omega^\bullet_\iota \to \Omega^\bullet_\kappa $
  such that $ \eta_{\kappa,\iota} (f) = H^* f$ for all $f\in \scrO_\iota (V) $
  with $V \subset U$ open. 
\end{lemma}

\begin{proof}
 By construction we have the following canonical identifications
\begin{equation}
\label{Eq:}
    \Omega^k_\iota (U) \cong  \Omega^k (\widetilde{U}) \big/ 
    \big( J_\iota \Omega^k (\widetilde{U}) +
    d J_\iota \wedge \Omega^{k-1} (\widetilde{U}) \big), \quad
    \text{where $ J_\iota :=  \scrI_\iota (\widetilde{U})$ }
\end{equation}
and
\begin{equation}
\label{Eq:}
    \Omega^k_\kappa (U) \cong  \Omega^k (\widehat{U}) \big/ 
    \big( J_\kappa \Omega^k (\widehat{U}) +
    d J_\kappa \wedge \Omega^{k-1} (\widehat{U}) \big), \quad
    \text{where $ J_\kappa :=  \scrJ_\kappa (\widehat{U})$}.
\end{equation}
Denote by $I \subset \calC^\infty (\widetilde{U})$ the vanishing ideal of $H(\widehat{U})$. The pull-back morphism 
$H^* : \Omega^k (\widetilde{U}) \to  \Omega^k (\widehat{U})$ then is surjective with kernel $I\Omega^k (\widetilde{U}) +   dI\wedge  \Omega^{k-1} (\widetilde{U})$, since $H(\widehat{U})$ is a closed submanifold of $\widetilde{U}$.
Since $H^*J_\iota  \subset J_\kappa$
and since $d$ commutes with $H^*$, pull-back by $H$ induces a surjective map
denoted by the same symbol
\[
 H^* : \Omega^k (\widetilde{U}) \big/ \big( J_\iota \Omega^k (\widetilde{U}) +
    d J_\iota \wedge \Omega^{k-1} (\widetilde{U}) \big) \longrightarrow 
    \Omega^k (\widehat{U}) \big/ \big( J_\kappa \Omega^k (\widehat{U}) +
    d J_\kappa \wedge \Omega^{k-1} (\widehat{U}) \big).
\]
This map is injective since
$I\Omega^k (\widetilde{U}) +   dI\wedge  \Omega^{k-1} (\widetilde{U})$ is contained in
$J_\iota \Omega^k (\widetilde{U}) + d J_\iota \wedge \Omega^{k-1} (\widetilde{U})$ and
since $H^*J_\iota = J_\kappa$. 
After choosing for each open set $V\subset U$ an open $\widehat{V} \subset \widehat{U}$ 
such that $V = \kappa^{-1}(\widehat{V})$, pull-back via $H_{|\widehat{V}}$ induces in
the same way isomorphisms $H_{V}^* : \Omega^k_\iota (V) \to \Omega^k_\kappa (V)$
for each $k\in \N$ and $V\subset U$ open. 
The family of isomorphisms $H_{V}^*$ then defines the desired sheaf isomorphism
$\eta_{\kappa,\iota}$. Since each $H_{V}^*$ commutes with the differentials,
$\eta_{\kappa,\iota}$ is a morphism of sheaf complexes, indeed. 
Moreover, since  $\Omega^\bullet_\iota$ is generated as a sheaf of differential graded 
algebras by $\Omega^0_\iota$,  $\eta_{\kappa,\iota}$ is uniquely determined by its action 
on $\Omega^0_\iota$. This finishes the proof.
\end{proof}

\begin{lemma}
\label{lem:step2}
Under the assumption of the preceding lemma let $G: \widecheck{U} \hookrightarrow \widetilde{U}$ be 
a second smooth embedding defined on an open neighborhood  $\widecheck{U} \subset \R^n$ 
of $\kappa (U)$ such that 
$\big( G_{|\kappa (U)},G^*\big) : (\kappa (U) , \scrO_\kappa ) \to  (\iota (U) , \scrO_\iota  ) $
is an isomorphism of locally ringed spaces . 
Then  $H_{V}^*  : \Omega^k_\iota (V) \to \Omega^k_\kappa (V)$ and $ G_{V}^* : \Omega^k_\iota (V) \to \Omega^k_\kappa (V)$
coincide for all open $V \subset U$ 
which means that the sheaf isomorphism $\eta_{\kappa,\iota}$ does not depend on the particular 
embedding inducing an isomorphism between $ (\kappa (U) , \scrO_\kappa )$ and $ (\iota (U) , \scrO_\iota  ) $.
\end{lemma}
\begin{proof}
After possibly shrinking $\widetilde{U}$ we can assume that $G (\widecheck{U})$ is closed in $\widetilde{U}$ as well.
For every smooth function $f \in \calC^\infty ( \widetilde{U}) $ we then have
\[
   (H^* f)_{|\widetilde{U} \cap \widecheck{U}} - (G^* f)_{|\widetilde{U} \cap \widecheck{U}}
   \in \scrJ_{\kappa} ( \widetilde{U} \cap \widecheck{U} ), \text{ and }
   (H^* d f)_{|\widetilde{U} \cap \widecheck{U}} - (G^* d f)_{|\widetilde{U} \cap \widecheck{U}} 
   \in d \scrJ_{\kappa} (\widetilde{U} \cap \widecheck{U}).
\]
That implies that the actions of $H^*$ and $G^*$ on $\Omega^k_\iota (U)$ coincide.
Likewise $H^*_{V} = G^*_{V}$ for all open $V\subset U$, hence $\eta_{\kappa,\iota}$
is independent of the particularly chosen embedding $H$.
\end{proof}

Next let  $\iota : U \hookrightarrow \R^n$
and $\kappa : V \hookrightarrow \R^m$ be two singular charts of $X$ defined on open 
$U,V \subset X$.  We will construct a sheaf morphism
$\eta_{\kappa,\iota}: \Omega^k_{\iota |U\cap V} \to \Omega^k_{\kappa |U\cap V}$. Let $x \in U\cap V$.
By the embedding theorem \ref{Thm:MinEmbeddingDim} there exists a singular chart
$\lambda : W_x \hookrightarrow \R^{\operatorname{rk} x}$ defined over an open neighborhood
$W_x\subset U\cap V$ of $x$. Moreover, after possibly shrinking $W_x$, there exist embeddings
$H: \widetilde{W}_x \hookrightarrow \R^n$
and $G: \widetilde{W}_x \hookrightarrow \R^m$ of an open neighborhood $\widetilde{W}_x$
of $\lambda (x)$ such that 
$\iota_{|W_x} = H\circ \lambda$ and $\kappa_{|W_x} = G\circ \lambda$, and such that 
$H^*$ induces an isomorphism from $(\iota (U) , \scrO_\iota )$ to $(\lambda (U) , \scrO_\lambda )$ and $G^*$ one
from $(\kappa (U) , \scrO_\kappa )$ to $(\lambda (U) , \scrO_\lambda )$. 
By Lemma \ref{lem:step1} we obtain isomorphisms of sheaf complexes
$\eta_{\lambda,\iota} : \Omega^\bullet_{\iota |W_x} \to \Omega^\bullet_{\lambda}$ and
$\eta_{\lambda,\kappa} : \Omega^\bullet_{\kappa |W_x} \to \Omega^\bullet_{\lambda}$.
Put $\eta_{\kappa,\iota}^{W_x} := (\eta_{\lambda,\kappa})^{-1} \circ \eta_{\lambda,\iota}$.
Then $\eta_{\kappa,\iota}^{W_x}$ is a sheaf isomorphism from
$\Omega^\bullet_{\iota |W_x}$ to $\Omega^\bullet_{\kappa |W_x}$ which by Lemma \ref{lem:step2} does not depend on the particular
choice of $\lambda$ and the embeddings $H$ and $G$. Moreover, if $y$ is another point of $U\cap V$,
an argument using Lemma \ref{lem:step2} and the embedding theorem \ref{Thm:MinEmbeddingDim} again shows that
the sheaf isomorphisms $\eta_{\kappa,\iota}^{W_x}$ and $\eta_{\kappa,\iota}^{W_y}$ coincide
on the overlap $W_x \cap W_y$. This proves the next lemma.

\begin{lemma}
 Given two singular charts $\iota : U \hookrightarrow \R^n$ and $\kappa : V \hookrightarrow \R^m$ of $(X,\scrO)$
 there exists a unique sheaf  morphism
 $\eta_{\kappa,\iota}: \Omega^\bullet_{\kappa |U\cap V} \to \Omega^\bullet_{\iota |U\cap V}$ such that 
 \[
   \eta_{\kappa,\iota |W_x} = (\eta_{\lambda,\kappa})^{-1} \circ \eta_{\lambda,\iota} 
 \]
 for each $x \in U\cap V$ and each singular chart $\lambda : W_x \hookrightarrow \R^{\operatorname{rk} x}$  defined on 
 a sufficiently small open neighborhood $W_x \subset U\cap V$ of $x$. 
\end{lemma}

Application of Lemma \ref{lem:step2} and the embedding theorem \ref{Thm:MinEmbeddingDim} 
a last time entails the final lemma.

\begin{lemma}
Assume that  $\iota : U \hookrightarrow \R^n$, $\kappa : V \hookrightarrow \R^m$, and $\lambda : W \hookrightarrow \R^l$ 
are three singular charts of $X$.   
Then the following cocycle condition holds true over the intersection $U\cap V \cap W$: 
\begin{equation}
  \label{eq:CocycleCondition}
  \eta_{\kappa , \iota} =
  \eta_{\kappa , \lambda} \circ   \eta_{\lambda,\iota} .
\end{equation}
\end{lemma}

The cocycle condition holding true for any triple of singular charts entails that the 
sheaf complexes $\Omega^k_\iota$ glue together to a global sheaf $\Omega^k_X$ 
on $X$. This sheaf is the  \emph{sheaf of abstract $k$-forms on} $X$. We 
sometimes denote it briefly by $\Omega^k$.
By construction, each of the $\eta_{\kappa,\iota}$ commutes with the exterior differential
$d$, hence the operators $d : \Omega^k_\kappa \to \Omega^{k+1}_\kappa$ glue together to a  sheaf morphism
$d : \Omega^k_X \to \Omega^{k+1}_X$. So finally we obtain a sheaf complex
$\big( \Omega^\bullet_X,d\big)$ of commutative differential   graded algebras.
The complex of global sections
$\big( \Omega^\bullet (X),d\big)$ will be called the \emph{Grauert--Grothendieck complex} of $X$.

\begin{remark}
  For $X \subset \C^n$ a complex space or algebraic variety, the complex $\Omega^\bullet (X)$ of abstract forms 
  on $X$ has been first constructed by Grauert \cite{GraKerDSKR} and Grothendieck \cite{GroRCAV}
  using sheaves of holomorphic respectively regular functions.    
  In \cite[\S 4.]{SpaDFDS}, Spallek has given a construction of the sheaf of $k$-forms of class $\scrC^\infty$ on an 
  affine differentiable space $(X,\scrO)$ which essentially corresponds to the one given here in 
  Eq.~\eqref{eq:def-diff-k-forms}. It is claimed in \cite{SpaDFDS} that given a (singular) atlas for 
  $(X,\scrO)$, the sheaves of $k$-forms over affine domains induce a global sheaf of $k$-forms on $X$. 
  Several details of the corresponding construction, in particular a verification of the cocycle condition, are missing. 
\end{remark}

\begin{remark}
 Recall that a morphism of differentiable spaces
 $(f,\varphi) : (X, \scrO_X) \to (Y, \scrO_Y)$ consists of
 a continuous map $f :X \to Y$ and a morphism of sheaves
 $\varphi : \scrO_Y \to f_* \scrO_X$. Using the construction of the 
 sheaves of abstract forms on $X$ and $Y$ one extends the morphism $\varphi$ 
 in  a  unique way to a morphism of sheaves of commutative 
 differential algebras
 $\varphi : \Omega^\bullet_Y \to f_* \Omega^\bullet_X$. One 
 concludes that forming the sheaf complex of abstract forms 
 is a functor defined on the category of differentiable spaces.
\end{remark}

%% file: basic-relative-forms.tex
\section{Basic relative forms}\label{sec:basic-relative-forms}
Let us first recall the notion of \emph{relative forms} associated to a smooth map $p: M \to N$ between manifolds $M$ and $N$.
By relative forms one understands sections of the sheaf complex 
$\Omega^\bullet_{M\overset{p}{\to} N}$ defined as the quotient sheaf 
\[
  \Omega^\bullet_M \big/ d (p^{-1}\scrC^\infty_N) \wedge \Omega^\bullet_M 
\] 
together with the differential induced by the differential on $\Omega^\bullet_M$. If $p:M\to N$ is a surjective submersion, 
the space of global sections $\Omega^k_{M\overset{p}{\to} N}(M)$ can be identified with the space of smooth families
$(\omega_y)_{y\in N}$ of forms $\omega_y \in \Omega^k \big( p^{-1} (y) \big)$. The differential acts fiberwise on 
$(\omega_y)_{y\in N}$ which means that 
\[ 
  d\big((\omega_y)_{y\in N} \big)= (d \omega_y)_{y\in N} \ .
\]
If the underlying map $p$ is projection onto the first factor of a product $N \times M$, one can identify 
$\Omega^k_{M\overset{p}{\to} N}$ with the sheaf of smooth sections $\Gamma^\infty ( -,\wedge^k s^* T^*M)$, where 
$s:N \times M \to M$ is projection onto the second factor. More precisely, one has in this case 
a sequence of sheaf morphisms, whose composition is an isomorphism: 
\begin{equation}
  \label{eq:two-sheaf-morphisms}
  \Gamma^\infty ( -,\wedge^k s^* T^*M) \longhookrightarrow \Omega^k_{N\times M} \longrightarrow 
  \Omega^k_{N\times M \overset{p}{\to} N} = \Omega^k_{N\times M} / d (p^{-1}\scrC^\infty_N) \wedge \Omega^{k-1}_{N\times M} \ . 
\end{equation}
Note that even though $\Gamma^\infty ( -,\wedge^k s^* T^*M)$ is a subsheaf of $\Omega^k_M$ for each $k$, one does not obtain
that way a subsheaf  complex of $\Omega^\bullet_M$ since the exterior derivative on $\Omega^\bullet_M$ 
does in general not map $\Gamma^\infty ( -,\wedge^k s^* T^*M)$ to $\Gamma^\infty ( -,\wedge^{k-1} s^* T^*M)$.
The correct differential on $\Gamma^\infty ( -,\wedge^k s^* T^*M)$ acts fiberwise as explained above.

After these preliminary remarks we now assume that $G$ is a compact Lie group acting on a smooth manifold $M$.
By $G\ltimes M$ we denote the corresponding action groupoid. 
More precisely, this is the Lie groupoid with arrow space $G\times M$,
object space $M$, source map $s: G\times M \to M$, $(g,p) \mapsto p$,
target map $t: G\times M \to M$, $(g,p) \mapsto gp$ and multiplication
$m : (G\ltimes M)^{(2)} \to G\times M$, $( (h,gp) , (g,p) ) \to (hg,p)$
defined on the fibered product $(G\ltimes M)^{(2)} := (G \times M) \sttimes (G\times M)$.  
The unit map of the action groupoid is $u: M \to G\times M$, $p \mapsto (e,p)$
with $e$ denoting the identity element of $G$, and the inversion map is 
$G\times M \to G\times M$, $(g,p) \mapsto (g^{-1},gp)$. 
We will denote the orbit space  of the  action groupoid $G\ltimes M$ by $M/G$, and the orbit map by 
$\pi : M \to M/G$. For more details on Lie groupoids see \cite{MoeMrcIFLG}.

The space $\inertianull (G\ltimes M)$ defined in Eq.~\eqref{eq:loop-space} corresponds 
to the loop space of the action groupoid  $G\ltimes M$ that means to the space 
of all $(g,p) \in  G\times M$ such that $ s(g,p) = t(g,p)$. 
Since $\inertianull (G\ltimes M)$ is a closed subspace of the manifold $G\times M$,
the loop space inherits from the ambient manifold the structure of a differentiable space. 
We denote the canonical embedding by $\iota : \inertianull (G\ltimes M) \hookrightarrow G \times M$.
Since by its definition the loop space is locally semialgebraic, it carries a
minimal Whitney B stratification, so $\inertianull (G\ltimes M)$ becomes a differentiable 
stratified space. For simplicity, we denote the loop space  shortly by 
$\inertianull$. Moreover, we denote the structure sheaf of smooth functions 
on $\inertianull$ by $\scrC^\infty_\inertianull$, and the sheaf of smooth functions 
on $G\times M$ vanishing on $\inertianull$ by $\scrJ$. 
Now consider the Grauert-Grothendieck complex $\Omega^\bullet_{\inertianull}$ of 
differential forms on the loop space. Following {\sc Brylinski}, we define the sheaf 
$\Omega^k_{\inertianull\to G}$ of \emph{relative forms on the loop space} as the quotient sheaf 
\begin{equation}
  \label{eq:definition-sheaf-relative-forms-loop-space}
  \Omega^k_{\inertianull\to G} :=  \iota^{-1} \Big( \Omega^k_{G \times M \to G} \big/ 
  \big(\scrJ\Omega^k_{G\times M \to G} + d\scrJ\wedge\Omega^{k-1}_{G\times M \to G}\big) \Big)\ .
\end{equation}
The graded sheaf $\Omega^\bullet_{\inertianull\to G}$ inherits from $\Omega^\bullet_{G\times M \to G}$ a differential 
turning it into a sheaf complex. Let us give another representation of the sheaf or relative forms on $\inertianull$.
To this end observe that the pull-back bundle $s^* T^*M$ induces a monomorphism of bundles 
$\wedge^k s^* T^*M \hookrightarrow \wedge^k T^* (G\times M)$, hence the  following morphism of sheaves:
\begin{equation}
\label{eq:sheaf-morphism}
  \Gamma^\infty_\inertianull ( - ,\wedge^k s^* T^*M)  \longrightarrow  \Omega^k_{\inertianull} = 
  \iota^{-1} \Big(\Omega^k_{G\times M} \big/ \big(\scrJ\Omega^k_{G\times M} + d\scrJ\wedge\Omega^{k-1}_{G\times M}\big)\Big).
\end{equation}
Here, $\Gamma^\infty_\inertianull (-,\wedge^k s^* T^*M)$ stands for the sheaf of smooth sections of the
vector bundle $\wedge^k s^* T^*M \to G\times M$ over the subspace $\inertianull$. 
In other words, the section space $\Gamma^\infty_\inertianull (U,E)$ for 
$U\subset \inertianull$ open can be identified with the quotient space 
$\Gamma^\infty(\widetilde{U},E)/ \scrJ(\widetilde{U}) \Gamma^\infty(\widetilde{U},E)$,
where $\widetilde{U} \subset G \times M$ open is chosen so that 
$\widetilde{U} \cap \inertianull = U$. 
Since the composition of sheaf morphisms in \eqref{eq:two-sheaf-morphisms} is an isomorphism, 
the sheaf morphism \eqref{eq:sheaf-morphism} induces a canoncial identification
\[
 \Omega^k_{\inertianull\to G} (U) \cong 
 \Gamma^\infty (\widetilde{U},\wedge^k s^* T^*M) \big/  
 \big(\scrJ(\widetilde{U})\ \Gamma^\infty (\widetilde{U},\wedge^k s^* T^*M) + 
 d\scrJ(\widetilde{U})\wedge \Gamma^\infty (\widetilde{U},\wedge^{k-1} s^* T^*M) \big),
\]
where $\widetilde{U}$ is chosen as before. For a section $\omega \in \Gamma^\infty (\widetilde{U},\wedge^k s^* T^*M)$
we denote  its image in $\Omega^k_{\inertianull\to G} (U)$ by $[\omega]$. 
Since $\inertianull$ is the union of the fibers 
$\{ g \} \times M^g$, $g\in G$, the relative form $[\omega]$ can be identified with the smooth family $(\omega_g)_{g\in G}$
of restrictions $\omega_g := \omega_{|M^g}$, and any smooth family  $(\omega_g)_{g\in G}$ of forms 
$\omega_g \in \Omega^k(M^g)$ gives rise to a unique relative form on the loop space. 
Under this identification the differential of $[\omega] = (\omega_g)_{g\in G}$  is given  
by the smooth family $(d \omega_g)_{g\in G}$.

Next recall that the $G$-action on $M$ gives rise for each $p\in M$ to the normal space $N_pM := T_pM/T_p\calO_p$,
where $\calO_p$ denotes the $G$-orbit through $p$. The family $N^*$ associating to each $p\in M$ the
conormal space $N_p^*M \subset T^*M$ is a smooth generalized vector subbundle of the cotangent bundle $T^*M$ 
in the sense of {\sc Drager--Lee--Park--Richardson} \cite{DraLeeParRicSDFG}. Note that under the
isomorphism between the tangent and cotangent bundle induced by a riemannian metric on $M$ 
the generalized vector bundle $N^*$  becomes a generalized distribution in the sense of 
{\sc Stefan--Suessmann}, cf.~\cite{SteASOFS,SusOFVFID}. The restriction of $N^*$ to an
orbit or a stratum of fixed orbit type is a vector bundle; see \cite{PflPosTanGOSPLG}. 
After these preparatory remarks it is clear that  $\wedge^kN^*$, $k\in \N^*$, is a smooth generalized 
vector subbundle of $\wedge^kT^*M$. Likewise, if one puts $N_pM^g := T_pM^g/(T_p\calO_p \cap T_pM^g)$ for $g\in G$, the 
alternating power $\wedge^k(NM^g)^*$ is a smooth generalized vector subbundle of
$\wedge^kT^*M^g$ for every $k\in \N^*$. The space $\Omega^k_{\textup{h.r.}} (U)$ of \emph{horizontal relative $k$-forms}
over $U\subset \inertianull$ open is now defined as 
\begin{equation*}
  \begin{split}
  \Omega^k_{\textup{h.r.}} (U) := \, & 
  \big\{ [\omega] = (\omega_g)_{g\in G} \in \Omega^k_{\inertianull\to G} (U) \mid \omega_{(g,p)} \in \wedge^kN^*_p
  \text{ for all $(g,p)\in \inertianull$} \big\} = \\ = &
  \big\{ [\omega] = (\omega_g)_{g\in G} \in \Omega^k_{\inertianull\to G} (U) \mid  \omega_g  \in 
  \Gamma^\infty (U\cap M^g , \wedge^k(NM^g)^* \big\} \ .
  \end{split}
\end{equation*}
Obviously, we thus obtain a subsheaf $\Omega^k_{\textup{h.r.}} \subset \Omega^k_{\inertianull\to G}$ having these 
spaces as its section spaces. 
Now observe that the $G$-action on the cotangent bundle $T^*M$ coming from the $G$-action on $M$  
induces a $G$-action on relative forms, hence we can speak of \emph{invariant relative $k$-forms}.
These are exactly those $[\omega] \in \Omega^k_{\textup{h.r.}} (U)$ which satisfy
\begin{equation}
  \label{eq:invariant-relativ-forms}
  \omega_{(hgh^{-1},hx)} ( hv_1, \ldots ,hv_k) =  \omega_{(g,x)} ( v_1, \ldots ,v_k)
\end{equation}
for all $(g,x) \in U$ and $h\in G$ such that $(hgh^{-1},hx) \in U$ and $v_1, \ldots ,v_k \in N_x$. Note that
by definition invariance of $[\omega]$ does not depend on the particular choice of a representative. If
one writes $[\omega]$ as a smooth family  $(\omega_g)_{g\in G}$ of forms $\omega_g$ on $U\cap M^g$ and if
$U$ is $G$-invariant, invariance of $[\omega]$ can be equivalently expressed by 
\begin{equation}
  \label{eq:invariant-relativ-forms-alternative}
  h^* \omega_{hgh^{-1}} =   \omega_g \quad \text{for all $g,h\in G$}.
\end{equation}
This relation implies in particular that the differential on  $\Omega^\bullet_{\inertianull\to G}$ maps
an invariant family  $(\omega_g)_{g\in G}$ over a $G$-invariant open $U$ to the invariant family 
$(d\omega_g)_{g\in G}$. Since the $G$-action on $T^*M$ leaves the conormal bundle $N^*$ invariant, 
we can even speak of \emph{invariant horizontal relative $k$-forms}. Now we are ready to put for 
$O \subset M/G$ open
\begin{equation}
  \label{eq:basic-relative-forms}
   \Omega^k_{\textup{b.r.}} (O) := \big\{ [\omega] \in \Omega^k_{ \textup{h.r.}} (s^{-1}_{|\inertianull}\pi^{-1}(O)) \mid 
   [\omega] \text{ is invariant} \big\}.
\end{equation}
These spaces are the section spaces of a sheaf $\Omega^k_{\textup{b.r.}}$. 
Observe that the sheaf $\Omega^k_{\textup{b.r.}}$ is defined over the orbit space $M/G$, not the loop space. 
Following {\sc Brylinski} again, we call sections of $\Omega^k_{\textup{b.r.}}$ \emph{basic relative $k$-forms}. In case we need 
to clarify the action groupoid underlying a sheaf of basic relative $k$-forms we will denote that sheaf more clearly by
$\Omega^k_{G\ltimes M \textup{-b.r.}}$.
The differential $d$ maps $\Omega^k_{\textup{b.r.}}$ to $\Omega^{k+1}_{\textup{b.r.}}$. This follows from Cartan's magic formula 
since it entails for $[\omega] = (\omega_g)_{g\in G} \in  \Omega^k_{\textup{b.r.}} (O) $,  $g\in G$,
and every element $\xi $ of the Lie algebra $\frakg_g$ of the centralizer $G_g := \mathsf{Z}_G (g)$ of $g$ 
the equality
\[
  i_{\xi_{M^{\!g}}} d \omega_g = \scrL_{\xi_{M^{\!g}}} \omega_g - di_{\xi_{M^{\!g}}} \omega_g  = 0,
\]
where $\xi_{M^{\!g}}$ denotes the fundamental vector field of $\xi$ on $M^g$. So we finally
obtain a complex of sheaves $(\Omega^\bullet_{\textup{b.r.}}, d)$ over the orbit space $M/G$.
Since each of the $\Omega^k_{\textup{b.r.}}$ is in a natural way a $\scrC^\infty_{M/G}$-module,
$(\Omega^\bullet_{\textup{b.r.}}, d)$ is even a complex of fine sheaves. To formulate our main result let us remind
the reader that we denote by  $\scrB$ Brylinski's sheaf over the orbit space $M/G$ and that this sheaf has 
section spaces given by \eqref{eq:brylinski-sheaf-section-spaces}. 

\begin{example}
  As an example let us consider the $\sphere^1$-action on $\R^2$ by rotation. We parametrize $\sphere^1$ by $e^{2\pi i \theta}$, 
  where $ \theta \in \R$. In coordinates, the action is expressed as 
  \[
    \R \times \R^2 \to \R^2, \: \big( \theta , (x,y) \big) \mapsto \big(\cos(\theta) x-\sin(\theta)y, \sin(\theta)x+\cos(\theta)y\big)\ . 
  \]
  The loop space is given by 
  \[
    \inertianull (\sphere^1\ltimes \mathbb{R}^2) = \big\{ (e^{2\pi i \theta}, (x,y)  ) \in \sphere^1\times \R^2 
    \mid \theta=0 \text{ or } x=y=0 \big\}.
  \]
  The minimal stratification of $\inertianull (\sphere^1\ltimes \mathbb{R}^2)$ is given by the decomposition into the strata 
  $S_0 = \{ 1 \} \times (\R^2  \setminus \{(0,0)\}$, $S_1 = ( \sphere^1 \setminus \{ 1 \} \times \{ (0,0) \} $, and 
  $S_2 = \{ (1,(0,0)) \} $.
  Using the given parametrization of $\sphere^1$ it is clear that in a neighborhood of the subset 
  $\{ 1 \} \times \R^2 \subset \inertianull (\sphere^1\ltimes \mathbb{R}^2)$  the loop space looks 
  like a neighborhood of $\{ 0 \} \times \R^2$ in the space 
  \[ 
   \Lambda_0' := \big( \R \times \{ (0,0)\} \big) \cup  \big( \{ 0 \} \times \R^2 \big) \ .
  \]  
  The loop space is smooth around each point of the stratum $S_1$, hence to describe the sheaf of smooth functions 
  on $\inertianull (\sphere^1\ltimes \mathbb{R}^2)$ we need to only understand how smooth funtions on $\Lambda_0'$ 
  look around a neighborhhod of the origin. To this end let $I \subset \scrC^\infty (\R^3)$ denote the ideal 
  of smooth functions vanishing on $\Lambda_0'$. If $f \in I$, then $f = \theta g$  for some $g \in \scrC^\infty (\R^3)$
  and with $\theta : \R^3 \to \R$ denoting here the projection onto the first coordinate. 
  Since $f$ vanishes on the $\theta$-axis, $g$ does so, too, hence $g = x f_1  + y f_2$ for some $f_1,f_2 \in \scrC^\infty (\R^3)$. 
  We obtain the representation $f = \theta x f_1 + \theta y f_2$. Therefore, the differential graded ideal 
  $I \Omega^\bullet  (\R^3) + dI \wedge  \Omega^\bullet (\R^3)\subset 
  \Omega^\bullet (\R^3) $ consists of all sums of forms of the form  
  \[
   \theta x  \omega_1 + \theta  y  \omega_2 + x  d\theta \wedge \omega_3 + y  d\theta \wedge \omega_4 + 
   \theta  dx \wedge \omega_3 + \theta dy \wedge \omega_4 , 
  \]
  where $ \omega_1, \omega_2   \in \Omega^k (\R^3)$ and $\omega_3, \omega_4  \in \Omega^{k-1} (\R^3)$. 
  Now observe that in the space of relative forms $\Omega^k_{\R^3 \overset{\theta}{\to}\, \R} (\R^3)$ the 
  form $d\theta$ vanishes. Moreover, the pullback of a relative form of degree $k \geq 1$ 
  to the stratum $\R^* \times \{ (0,0\}$ vanishes as well. 
  One concludes that, for $k \geq 1$ the space 
  \[
    \Omega^k_{\Lambda_0' \overset{\theta}{\to}\, \R} (\Lambda_0') := 
    \Omega^k_{\R^3 \overset{\theta}{\to}\, \R} (\R^3) \big/ I \Omega^k_{\R^3 \overset{\theta}{\to}\, \R}  (\R^3) + 
    dI \wedge  \Omega^{k-1}_{\R^3 \overset{\theta}{\to}\, \R} (\R^3)
  \] 
  of relative $k$-forms on $\Lambda_0'$ can be identified with the space of smooth families 
  $(\omega_\theta)_{\theta \in \R}$, where 
  \[
    \omega_\theta \in
    \begin{cases}
      \Omega^k ( \R^2 ), & \text{if $\theta =0$}, \\
      \{ 0 \}, & \text{else}.
    \end{cases}
  \]
  So one concludes that  
  $\Omega^k_{\inertianull (\sphere^1\ltimes \mathbb{R}^2) \to \sphere^1}\big(\inertianull (\sphere^1\ltimes\mathbb{R}^2)\big)
  \cong \Omega^k (\R^2)$ for $k \geq 1$ and that, under this ismorphism,
  $\Omega^k_{\textup{b.r.}} \big( \R^2/\sphere^1 \big) \cong \Omega^k_{\textup{bas}}\big( \R^2/\sphere^1 \big)$,
  where the latter denotes the space of basic $k$-forms on $\R^2$.
  In case $k=0$, one has 
  $\Omega^0_{\textup{b.r.}} \big(\R^2/\sphere^1 \big) \cong \scrC^\infty (\R^2/\sphere^1)$.
  Finally in this example, Brylinski's sheaf can be identified with the sheaf (over the orbit space) 
  of $G$-invariant locally constant functions on $\R^2$.
\end{example} 

\begin{theorem}
\label{thm:ResBrylinskiSheaf}
 Let $G$ be a compact Lie group acting on a manifold $M$. 
 The sheaf complex $(\Omega^\bullet_{\textup{b.r.}}, d)$ of basic relative forms 
 together with the natural monomorphism of sheaves 
 \begin{equation*}
   \label{eq:Brylinski-resolution}
    d_{-1}: \scrB \hookrightarrow \Omega^0_{\textup{b.r.}} = \pi_* \big(\scrC^\infty_\inertianull \big)^G 
 \end{equation*}
 then forms a fine resolution of  Brylinski's sheaf $\scrB$. 
\end{theorem}
\begin{proof}
  Let $O\subset M/G$ be open, and $f \in \scrB (O)$. By definition through Eq.~\eqref{eq:brylinski-sheaf-section-spaces}
  $f$ then is a smooth $G$-invariant function on $s^{-1}_{|\inertianull} \pi^{-1} (O)$, and $f(g,-):M^g \to\R$ is
  locally constant for all $g\in G$. Hence $df(g,-) \in \Omega^1 (M^g)$ vanishes for every $g$. This entails that 
  $\scrB \hookrightarrow \Omega^\bullet_{\textup{b.r.}}$ is a cochain complex of sheaves over the orbit space $M/G$. 
 
  It remains to show that that cochain complex of sheaves is exact, meaning that 
  for each orbit $\scrO \in M/G$ the complex of stalks $\scrB_\scrO \hookrightarrow \Omega^\bullet_{\textup{b.r.},\scrO}$ is exact.
  To this end we proceed in several steps. 
  In the first step we consider the case, where $M$ is a finite dimensional vector space $V$ carrying a linear $G$-action,
  and where $\scrO = \{0 \}$, the orbit through the origin. Choose a $G$-invariant scalar product on $V$. Assume that $B \subset V$ 
  is an open ball around the origin, and consider the homothety $h: [0,1] \times V \to V$, $(t,v) \to tv$. The homothety $h$ leaves 
  each of the subsets $B^g \subset B$ invariant, and commutes with the $G$-action. Hence 
  \[
   K : \Omega^k_{\textup{b.r.}} (B) \to \Omega^{k-1}_{\textup{b.r.}} (B) ,\:
   \omega = (\omega_g)_{g\in G} \mapsto K \omega :=
   \begin{cases}
     \inertianull (G \ltimes M) \ni (g,p) \mapsto \omega_g (0), & \text{for $k=0$},\\
     \big( \int_0^1 h_t^* (\xi_t \amslrcorner \omega_g) \big)_{g\in G} , & \text{for $k\geq 1$},
   \end{cases}
  \]
  is a well-defined operator, where $\Omega^{-1}_{\textup{b.r.}}$ is Brylinski's sheaf $\scrB$, $h_t$ equals $h (t,-):  B \to B$, $v \mapsto tv$, 
  and $\xi_t : B \to TB$ is the vector field given by $\xi_t := \partial_th_t$. Cartan's magic formula implies that
  \[
     (\omega_g)_{g\in G} = d K (\omega_g)_{g\in G}  + Kd (\omega_g)_{g\in G} 
     \quad \text{for all $(\omega_g)_{g\in G} \in \Omega^k_{\textup{b.r.}} (B)$ and $k\in \N$},
  \]
  since $h_1^*$ acts by identity on $\big( \omega_g \big)_{g\in G}$, $\big( h_0^* \omega_g \big)_{g\in G} =0$ for $k\geq 1$, and 
  $h_0^* \omega = \omega (-,0)$ for $k=0$. Hence $\scrB_\scrO \hookrightarrow \Omega^\bullet_{\textup{b.r.},\scrO}$ is exact when $\scrO = \{0\}$.

  In the second step we come back to the general case of a $G$-action on an arbitrary manifold $M$.
  Choose a $G$-invariant riemannian metric $\eta$ on $M$. Let $p\in M$ be a point, $\scrO$ the orbit through $p$, and $N_p := T_pM/T_p\scrO$
  the normal space to $\scrO$ at $p$. Via the riemannian metric we can identify $N_p$ with the orthogonal complement of $T_p\scrO$ in $T_pM$. 
  The isotropy group $G_p$ acts in a natural way on $N_p$, cf.~\cite[Sec.~4.2.5]{PflAGSSS}.
  Choose an open ball $B$ around the origin of $N_p$ with radius smaller than the injectivity radius at $p$.  
  The slice theorem \cite[Thm.~4.2.6]{PflAGSSS} entails that the $G$-action on $M$ induces an action of the isotropy group $G_p$ on the 
  slice $S_p := \exp (B)$ and that the exponential map intertwines the $G_p$-actions 
  on $N_p$ and $S_p$. Moreover, the exponential map procvides an equivariant diffeomorphism between $G \times_{G_p} B$ 
  and the $G$-saturation $U : = G\cdot S_p$ of the slice. It therefore suffices to verify the claim for the case where $M$ has the form 
  $G\times_H B$ with $H \subset G$ being a compact subgroup and $B$ an open $H$-invariant ball around  the origin of a finite dimensional
  $H$-representation space $V$, and where $\scrO$ is the orbit through the point $[e,0] \in G\times_H B$.
  From now on we will consider only this setting. Note that here and in the following we will denote by $[g,v]$ the equivalence class
  of a point $(g,v) \in G\times B$ in $G\times_H B$.

  In the third step we provide a description of the tangent bundle $T (G\times_H V)$. To this end choose a bi-invariant riemannian metric 
  on $G$. Let $\scrQ$ be the foliation of $G$ given by the orbits of the canonical right action of $H$ on $G$. 
  For each $g\in G$ let $E_g$  be the orthogonal complement of the tangent space $T_g\scrO$ to the 
  leaf of $\scrQ$ through $g$. One thus obtains a vector bundle $E\to G$ which is invariant under the left action of $G$ on $TG$ 
  and invariant under the right action of $H$. The latter follows from the fact that the right action of $H$ on $G$ maps leaves of $\scrQ$ 
  to leaves, since $g \exp (t\xi) h = gh \exp\big( t \operatorname{Ad}_{h^{-1}} (\xi)\big)$ for all $g\in G$, $h\in H$, $\xi \in \frakh$, and $t\in \R$. 
  One concludes that $E$ can be identified with the trivial bundle $G \times \frakm$, where $\frakm$ is the orthogonal complement of
  $\frakh$ in $\frakg$. Now call two elements $(\Xi,(v,X)), (\Zeta ,(w,Y)) \in E \times TV = E \times V \times V$ 
  \emph{equivalent}, if there is an $h\in H$ with $\Xi = \Zeta h $ and $h(v,X) = (w,Y)$.
  We denote by $[\Xi,(v,X)]$ the equivalence class of  $(\Xi,(v,X)) \in E \times TV$.
  One checks immediately that the quotient space of  $E \times TV$ by this equivalence relation can be canonically 
  identified  with the tangent bundle $T (G\times_H V)$. Under this identification, an element of the tangent space 
  $T_{[g,v]} (G\times_H V)\cong \frakm \times V$ over the footpoint $[g,v] \in G\times_H V$ has  a unique representation of 
  the form $[(g ,\xi),(v,X)] $ with $\xi \in \frakm$ and $X \in V$. For later purposes let us remark that 
  if $[g,v] = [g',v']$, i.e.~if $(g',v') = (gh^{-1},hv) $ for some $h\in H$, then 
  $[( g,\xi),(v,X)] = [( g',\xi'),(v',X')] $ with $\xi' = \xi h^{-1}$ and $v' = hv$. In the following step
  we will denote the equivalence class $[( g,\xi),(v,X)]$ shortly by $[\xi,X]_{(g,v)}$. 
 
  The fourth step consists in verifying that the embedding  map $\iota: B \hookrightarrow M:= G\times_H B$, $v \mapsto [e,v]$ 
  induces an isomorphism between the sheaves $\Omega^k_{G \ltimes M \textup{-b.r.}}$ and $\Omega^k_{H \ltimes B \textup{-b.r.}}$. Observe that 
  both sheaves live on the same topological space, since the orbit spaces $M/G$ and $B/H$ are canonically isomorphic since $\iota$ is a 
  Morita equivalence.  The isomorphism is given by pullback via $\iota$. More precisely, for $O \subset M/G$ 
  and a basic relative form $\omega =(\omega_g)_{g\in G} \in \Omega^k_{G\ltimes M\textup{-b.r.}} (O)$ let  $\iota^*\omega$
  be the family $\big(\iota^*_{O,h} \omega_h\big)_{h\in H}$, where
  $\iota_{O,h}$ denotes the restriction of the embedding $\iota$ to $\big( \pi^{-1}_{H\ltimes B} (O) \big)^h$ with $\pi_{H\ltimes B} : B \to B/H$
  being the orbit map of the groupoid $H\ltimes B$. Obviously, $\iota^*_{O,h} \omega_h$ is a $\mathsf{Z}_H (h)$-invariant horizontal form on
  $\big( \pi^{-1}_{H\ltimes B} (O) \big)^h $, and the family $\big(\iota^*_{O,h} \omega_h\big)_{h\in H}$ is $H$-invariant. Hence  
  $\iota^*\omega \in \Omega^k_{H\ltimes B\textup{-b.r.}} (O)$, so we obtain a morphism of sheaf complexes 
  $\iota^*:\Omega^\bullet_{G\ltimes M\textup{-b.r.}} \to \Omega^\bullet_{H\ltimes B\textup{-b.r.}}$. Let us show that it is an isomorphism. 
  To this end note first that for every $g\in G$ the invariant space $M^g$ is 
  given by $M^g = \{ [f,v] \in G\times_HB\mid g \in fH_vf^{-1}\}$. 
  Now let $(\varrho_h)_{h\in H} \in \Omega^k_{H \ltimes B \textup{-b.r.}}(O)$, and 
  define $\omega \in \Omega^k_{G \ltimes M \textup{-b.r.}}(O) $ by 
  \[
    \omega_{g,[f,v]}\big( [\xi_1,X_1]_{(f,v)}, \ldots, [\xi_k,X_k]_{(f,v)} \big) := 
    \varrho_{f^{-1}gf ,v} (f^{-1}gfX_1 ,\ldots, f^{-1}gf X_k ) ,
  \]
 where $g\in G$, $[f,v]\in M^g \cap s_{\inertianull (G \ltimes M)}^{-1}\pi_{|G \ltimes M}^{-1} (O)$, and 
 $[\xi_1,X_1]_{(f,v)}, \ldots, [\xi_k,X_k]_{(f,v)} \in T_{(f,v)}  M$.
 Since $\varrho = (\varrho_h)_{h\in H}$ is  $H$-invariant, one obtains for 
 every $h\in H$ the equality
 \[ 
  \begin{split}
   \omega_{g,[fh^{-1},hv]} \, & \big( [\xi_1h^{-1},hX_1]_{(fh^{-1},hv)}, \ldots, 
   [\xi_kh^{-1},hX_k]_{(fh^{-1},hv)} \big)    =  \\
   = \, & \varrho_{hf^{-1}gfh^{-1} ,hv} (hf^{-1}gfX_1 ,\ldots, hf^{-1}gf X_k ) = \\
   = \, & \varrho_{f^{-1}gf ,v} (f^{-1}gfX_1 ,\ldots, f^{-1}gf X_k ).
  \end{split}
 \] 
 This shows that $\omega_g$ is independant of the choices made, and an element of 
 $\Omega^k (M^g \cap s_{\inertianull (G \ltimes M)}^{-1}\pi_{|G \ltimes M}^{-1} (O))$. 
 Moreover, $\omega_g$ is a horizontal form by construction. The family
 $\omega = (\omega_g)_{g\in  G}$ is also $G$-invariant. To verify this let $h\in G$, and 
 observe that by definition
 \[
    \omega_{hgh^{-1},[hf,v]}\big( [h\xi_1,X_1]_{(hf,v)}, \ldots, [h\xi_k,X_k]_{(hf,v)} \big) = 
    \varrho_{f^{-1}gf ,v} (f^{-1}gfX_1 ,\ldots, f^{-1}gf X_k ),
  \]
  where $g\in G$, $[f,v]$, and $[\xi_1,X_1]_{(f,v)}, \ldots, [\xi_k,X_k]_{(f,v)}$ are as above.
  This proves $G$-invariance of the family $\omega$, hence $\omega \in \Omega^\bullet_{G\ltimes M\textup{-b.r.}} (O)$ indeed. 
  By construction it is clear that $\iota^* \omega = \varrho$. By $G$-invariance, $\omega$ is uniquely determined by
  $\iota^* \omega$.  So $\iota$ is a sheaf isomorphism as claimed.   
  
  In the fifth and final step we show that $\scrB_\scrO \hookrightarrow \Omega^\bullet_{G\ltimes M\textup{-b.r.},\scrO}$ is an exact sheaf complex
  in the case where $M = G\times_H B$ and $\scrO$ is the orbit through the point $[e,0]$. By the second step, it suffices to
  consider this case. By the fourth step, the embedding $\iota : B \hookrightarrow M$ induces an isomorphis of sheaf complexes
  $\iota^* :  \Omega^\bullet_{G\ltimes M\textup{-b.r.}}\to \Omega^\bullet_{H\ltimes B\textup{-b.r.}} $. By the first step, 
  the cochain complex $\scrB_{\{0\}} \hookrightarrow \Omega^\bullet_{H\ltimes B\textup{-b.r.},\{ 0 \}}$ is exact, hence 
  $\scrB_\scrO \hookrightarrow \Omega^\bullet_{G\ltimes M\textup{-b.r.},\scrO}$ is so, too, since the orbit through $\iota (0)$ coincides with
  $\scrO$. The proof is finished. 
\end{proof}


%% file: appendix.tex
\section{Differentiable stratified spaces}
For the convenience of the reader we briefly recall here the notion of a 
differentiable space, mainly following \cite{NGonzalezSanchoBook}, 
and then describe what it means that a stratification is compatible with the 
differentiable structure.

\begin{definition}
  An algebra over $\R$ of the form $A =\scrC^\infty (\R^n)/ J$, where 
  $J \subset \scrC^\infty (\R^n)$ is a closed ideal, is called a
  \emph{differentiable algebra}. 
  By $\spec A$ is the maximal spectrum of a differentiable algebra $A$,
  and by $\scrO_A$ its structure sheaf, which is the sheafification
  of the presheaf $U \mapsto A_U$, $U\subset \spec A$ open, where 
  $A_U$ is the localization of $A$ over the subset of elements 
  which do not vanish over $U$; cf.~\cite[Sec.~3.1]{NGonzalezSanchoBook}.  
\end{definition}
  A differentiable algebra carries in a natural way the structure of a Fr\'echet algebra.
  Moreover, given a differentiable algebra $A$, the pair $(\spec A, \scrO_A)$  
  is a commutative locally ringed space. 
\begin{definition}
  A commutative locally $\R$-ringed space $(X,\scrO)$ is called an 
  \emph{affine differentiable space}, if it is isomorphic as a 
  commutative locally ringed space to $(\spec A, \scrO_A)$ 
  for some differentiable algebra $A$. 
  The ringed space $(X,\scrO)$ is called a \emph{differentiable space}, if 
  for every $x\in X$ there exists an open neighborhood $U \subset X$ such that 
  $(U,\scrO_{|U})$ is an affine differentiable space. By a \emph{morphism of differentiable spaces}
  we understand a morphism of locally $\R$-ringed spaces
  $(f,\varphi) : (X, \scrO_X) \to (Y, \scrO_Y)$ 
  between differentiable spaces $(X, \scrO_X)$ and $(Y, \scrO_Y)$. 
\end{definition}
 
  Locally, a differentiable space $(X,\scrO)$ can be embedded into euclidean space. 
  We call an embedding $\iota : U \hookrightarrow \R^n$ over an open $U \subset X$ 
  together with a morphism of differentiable algebras $\iota^* : \scrC^\infty (\iota (U)) \to \scrO (U)$
  such that $\big( \iota^*  (x_1) (x) ,  \ldots , \iota^*  (x_n) (x) \big)=  \iota (x)$
  for all $x \in U$ a \emph{singular chart} for $(X,\scrO)$ if these data induce an 
  isomorphism of locally $\R$-ringed spaces 
  \[
    (\iota,\iota^*) : (U,\scrO_{|U}) \rightarrow \big( \iota (U), \scrC^\infty_{|\iota (U)}/\scrJ_\iota \big), 
  \]
  where $\scrJ_\iota$ is the kernel of the morphism of sheaves  
  $\iota^* : \scrC^\infty_{|\iota (U)} \to \scrO_{|U}$.
  By a \emph{singular atlas} of  $(X,\scrO)$ we understand a 
  family $\mathfrak{A}$ of singular charts  such that the family of domains 
  $\{ \operatorname{dom} (\iota) \}_{(\iota,\iota^*)\in \mathfrak{A}}$ is an open cover of $X$.
  If  $(X,\scrO)$ is reduced, the embedding $\iota : U \hookrightarrow \R^n$ completely determines 
  the singular chart $(\iota,\iota^*)$, because $\iota^*$ then is just pullback by the embedding. 
  By abuse of language we  denote a singular chart even in the non-reduced case just by the embedding
  $\iota : U \hookrightarrow \R^n$.

  The following result shows that around a point of a differentible space the minimal embedding dimension is given
  by the dimension of the Zariski tangent space. 

\begin{theorem}[{\cite[Prop.~1.3.10 \& Corollaries]{PflAGSSS}}]
\label{Thm:MinEmbeddingDim}
  Let $(X,\scrO)$ be differentiable space, and  $x\in X$  a point.
  Then there exists an open affine neighborhood $W$ of $x$ together with a singular chart
  $\lambda : W \hookrightarrow \R^{\operatorname{rk} x}$, where $\operatorname{rk} x$ is the dimension of 
  the Zariski tangent space $T_xX$ at $x$. Moreover, if $\iota : U \hookrightarrow \R^n$ is
  another singular chart defined on an open affine neighborhood of $x$, then $\operatorname{rk} x \leq n$, and
  there exists an open affine neighborhood $V \subset W\cap U$ of $x$, an open neighborhood
  $\widetilde{V} \subset \R^{\operatorname{rk} x}$ of $\lambda (x)$
  and a smooth embedding $H : \widetilde{V} \hookrightarrow \R^n$ such that
  $\iota (V)$ is closed in  $\widetilde{V}$ and such that $H\circ \iota_{|V} =\kappa_{|V}$.
\end{theorem}

\begin{definition}
  A stratification $\scrS$ of the topological space $X$ underlying a differentiable space $(X,\scrO)$ is said to be
  \emph{compatible with} $(X,\scrO)$ or just \emph{with} $\scrO$ if for each stratum $S\in \scrS$ and singular chart 
  $\iota: U \hookrightarrow \R^n$ the image $\iota (S\cap U)$ of the stratum under $\iota$ is a submanifold of $\R^n$. 
  We call a differentiable space $(X,\scrO)$ together with a compatible stratification $\scrS$ of $X$ 
  a \emph{differentiable stratified space}.
\end{definition}

\begin{example}
  Typical examples of differentiable stratified spaces are real or complex algebraic varieties and orbit spaces 
  of compact Lie group actions on manifolds. In \cite{PflPosTanGOSPLG} it has been shown that the orbit space of a proper Lie groupoid
  carries the structure of a differentiable stratified space in a natural way.
\end{example}
